\documentclass[12pt]{article}
\usepackage{graphicx}
\usepackage{amsmath}
\usepackage[dvips]{epsfig}
\usepackage{amssymb}

\makeatletter
\renewcommand{\section}{\@startsection
  {section}%
  {2}%
  {0mm}%
  {\baselineskip}%
  {0.5 \baselineskip}%
  {\centering}}
\makeatother
\begin{document}

\title { On  $q$-Bernstein and $q$-Hermite polynomials}
\author{  T. Kim$^1$, J. Choi$^1$, Y. H. Kim$^1$, C. S.  Ryoo$^2$  \\[0.5cm]
$^1$Division of General Education-Mathematics,\\
 Kwangwoon University, Seoul 139-701,  Korea\\
 \\
           $^2$Department of Mathematics, \\
          Hannam University, Daejeon 306-791, Korea
   }

\date{}
\maketitle

 {\footnotesize {\bf Abstract}\hspace{1mm}
Recently, T. Kim([1,2]) introduced  $q$-Bernstein  polynomials which
are different $q$-Bernstein  polynomials of Phillips.  In this
paper, we investigate some identities for the $q$-Bernstein
polynomials which are treated by Kim in previous paper([1]).



\bigskip
\section{Introduction }
\bigskip
Let $C[0,1]$ denote the set of continuous functions on $[0, 1]$.
For $f \in C[0,1]$,  Bernstein introduced the following well-known
linear operator(see [3, 5, 8]):
$$ \aligned  \Bbb B_{n}(f \mid x) &=\sum_{k=0}^n  f \left( \dfrac{k}{n}\right)
\binom{n}{k} x^k (1-x)^{n-k} \\
&=\sum_{k=0}^n f \left( \dfrac{k}{n}\right)B_{k,n}(x),
\endaligned  \eqno(1)$$
where $ \binom{n}{k}= \dfrac{n(n-1) \cdots (n-k+1)}{n!}$.  Here $
\Bbb B_{n}(f \mid x)$ is called the Bernstein operator of order
$n$ for $f$. For $n, k \in \Bbb Z_+ = \Bbb N  \cup \{0\},$
 the Bernstein polynomial of degree $n$ is defined by
$$ B_{k, n}(x)=\binom{n}{k} x^k (1-x)^{n-k}, \text{ (see [1-8]) }. $$

The Hermite polynomials have the generating function
$$ e^{2zw-w^2}= \sum_{n=0}^\infty \dfrac{H_n(z)}{n!} w^n, \quad z, w \in \Bbb C, \eqno(2) $$
which gives the Cauchy type integral
$$H_n(z)= \dfrac{n!}{ 2 \pi i } \int_C e^{2zw-w^2} \dfrac{d^w}{w^{n+1}}, $$
where $C$ is a circle around the origin with the
positive direction (see [4]).

From (2), we note that the Hermite polynomials are given by
$$H_n(x)=n! \sum_{k=0}^{\left[ \frac{n}{2}\right]}  \dfrac{(-1)^k}{k! (n-2k)!} (2x)^{n-2k}, $$
where [ $\cdot$ ] is the Gauss symbol (see [4]).

Let $ 0< q<1$. Define the $q$-number of $x$ by
$$[x]_q =\frac{1-q^x}{1-q}.  $$
Note that $\lim_{q \rightarrow 1} [x]_q=x$ (see [9-11]).

For $ f \in C[0,1]$, the $q$-Bernstein type  operator of order $n$ for
$f$ is defined by
$$ \aligned   \Bbb B_{n, q}(f \mid x) & =\sum_{k=0}^n  f \left( \dfrac{k}{n}\right) \binom{n}{k} [x]_q^k
[1-x]_{\frac{1}{q}}^{n-k}\\
&=\sum_{k=0}^n \left( \dfrac{k}{n}\right) B_{k,n}(x \mid q),
\endaligned \eqno(4)
$$ where $n, k \in \Bbb Z_+$ (see [1]). Here  $ B_{k,n}(x \mid q)= \binom{n}{k}[x]_q^k
[1-x]_{\frac{1}{q}}^{n-k}$ are called the $q$-Bernstein type
polynomials of degree $n$ (see [1, 2]). Note that (4) is
a $q$-analogue of (1).

In this paper, we consider the $q$-Hermite polynomials related to the
$q$-Bernstein type polynomials. From these $q$-Hermite
polynomials, we give some interesting identities for the
$q$-Bernstein polynomials.

\bigskip
\section{$q$-Bernstein polynomials and  $q$-Hermite polynomials }
\bigskip

In this section, we consider a $q$-analogue of (2), which are called
the $q$-Hermite polynomials, as follows:
$$ e^{2t[x]_q-t^2}= \sum_{n=0}^\infty H_{n, q}(x) \dfrac{t^n}{n!}. \eqno(5) $$
From (5), we have

$$ \aligned   e^{2t[x]_q-t^2} &=  \left(  \sum_{m=0}^\infty  \dfrac{(2[x]_q)^m}{m!} t^m \right)
\left(   \sum_{l=0}^\infty   \dfrac{(-1)^l}{l!} t^{2l}\right)\\
&=\sum_{n=0}^\infty  \left(   \sum_{l=0}^{\left[
\frac{n}{2}\right]}  \dfrac{(-1)^l(2[x]_q)^{n-2l} n!}{l! (n-2l)!}
 \right) \dfrac{t^n}{n!} .\endaligned \eqno(6)
$$
Therefore, by (5) and (6), we obtain the following theorem.

\bigskip

{ \bf Theorem 1.}  For $ n \in \Bbb Z_+$, we have
$$ H_{n,q}( x )=  n! \sum_{l=0}^{\left[
\frac{n}{2}\right]}  \dfrac{(-1)^l 2^{n-2l} [x]_q^{n-2l}}{l!
(n-2l)!}.$$

\bigskip
Note that Theorem 1 is a $q$-analogue of (3).  By (4), we obtain the
following corollary.

\bigskip

{ \bf Corollary 2.}  For $ x \in [0, 1]$, we have
$$ H_{n,q}( x )=  n! \sum_{l=0}^{\left[
\frac{n}{2}\right]}  \dfrac{(-1)^l 2^{n-2l} B_{2l, n}(1-x \mid
\frac{1}{q})}{l! (n-2l)! [1-x]_{\frac{1}{q}}^{2l}}.$$

\bigskip

From the definition of $q$-Bernstein polynomials, we note that
$$ \aligned    \sum_{n=0}^\infty  \dfrac{(-1)^n}{n!} \left( \dfrac{d}{d[x]_q} \right)^n e^{-t [x]_q}
&=   \sum_{n=0}^\infty  \dfrac{(-1)^n}{n!} (-t)^n e^{-t [x]_q}\\
&= \left(\sum_{n=0}^\infty  \dfrac{t^n}{n!}  \right) e^{-t [x]_q}= e^{[1-x]_{\frac{1}{q}} t} \\
&= \dfrac{k!}{ (t[x]_q)^k} \left( \dfrac{(t[x]_q)^k e^{[1-x]_{\frac{1}{q}} t} }{k!}\right)\\
&= \dfrac{k!}{ (t[x]_q)^k} \sum_{n=k}^\infty B_{k, n}(x \mid q)\dfrac{t^n}{n!} \\
&= \dfrac{k!}{ (t[x]_q)^k} \sum_{n=0}^\infty B_{k, n+k}(x \mid q)\dfrac{t^{n+k}}{(n+k)!} \\
&=\dfrac{1}{ [x]_q^k} \sum_{n=0}^\infty \dfrac{k!}{(n+k)\cdots
(n+1)} B_{k, n+k}(x \mid
q)\dfrac{t^n}{n!} \\
&= \dfrac{1}{ [x]_q^k} \sum_{n=0}^\infty \dfrac{B_{k, n+k}(x \mid
q)}{\binom {n+k}{k}} \dfrac{t^n}{n!} .\endaligned \eqno(7)
$$
By a simple calculation, we easily get
$$ \aligned    \sum_{n=0}^\infty  \dfrac{(-1)^n}{n!} \left( \dfrac{d}{d[x]_q} \right)^n e^{-t [x]_q}
&=   \sum_{n=0}^\infty  \left(  \sum_{l=0}^n \dfrac{(-1)^l n!  }{l!(n-l)!}\left( \dfrac{d}{d[x]_q} \right)^l [x]_q^{n-l} \right)\dfrac{t^n}{n!} \\
&=  \sum_{n=0}^\infty \left( \sum_{l=0}^n \binom nl (-1)^l \left(
\dfrac{d}{d[x]_q} \right)^l [x]_q^{n-l} \right)\dfrac{t^n}{n!}
.\endaligned \eqno(8)
$$
Therefore, by (7) and (8), we obtain the following theorem.

\bigskip

{ \bf Theorem 3.}  For $  x \in [0, 1], n, k  \in \Bbb Z_+$, we
have
$$ [x]_q^k \sum_{l=0}^n \binom nl (-1)^l \left(
\dfrac{d}{d[x]_q} \right)^l [x]_q^{n-l}=  \dfrac{B_{k, n+k}(x \mid
q)}{\binom {n+k}{k}} .$$

\bigskip

By a simple calculation, we easily get
$$ \left( \dfrac{1}{2} \dfrac{d}{d[x]_q} \right)^n e^{-2t [x]_q} =(-t)^n e^{-2t[x]_q}.$$
Thus, we have
$$\sum_{n=0}^\infty  \dfrac{(-1)^n}{n!}\left( \dfrac{1}{2} \dfrac{d}{d[x]_q} \right)^{2n}
 e^{-2t [x]_q}=   e^{-t^{2n}+2t[x]_q}
 e^{-4t [x]_q}. \eqno(9)$$

From (6) and (9), we can derive the following equation:

$$  e^{-t^{2n}-2t[x]_q}= \sum_{n=0}^\infty \left( n! (-1)^n   \sum_{l=0}^{\left[
\frac{n}{2}\right]}  \dfrac{(-1)^l 4^{n-l} [x]_q^{n-l} }{l!
(n-2l)!} \right) \dfrac{t^n}{n!} \eqno(10) $$ and
$$  e^{-t^{2}+2t[x]_q}e^{-4t [x]_q}= \sum_{n=0}^\infty \left(   \sum_{l=0}^{n}
\binom nl H_{l,q}(x) (-1)^{n-l} [x]_q^{n-l}\right)
\dfrac{t^n}{n!}. \eqno(11) $$

Therefore, by (10) and (11), we obtain the following theorem.

\bigskip

{ \bf Theorem 4.}  For $   n  \in \Bbb Z_+$, we have
$$ H_{n,q}(x)= \sum_{l=0}^n   \binom nl H_{l,q}(x) (-1)^{n-l} [x]_q^{n-l}.$$

\bigskip

 By (4) and Theorem 4, we obtain the following corollary.

\bigskip

{ \bf Corollary 5.}  For $   n  \in \Bbb Z_+$  and $  x \in [0,
1]$, we have
$$ H_{n,q}(x)= \sum_{l=0}^n   \binom nl H_{l,q}(x) (-1)^{l} B_{l,n}(1-x \mid \frac{1}{q}) [1-x]_{\frac{1}{q}}^{-l}.$$

\bigskip

For $ t \in \Bbb R$ and $  x \in [0, 1]$, we have
$$ \aligned   \sinh ( t[1-x]_{\frac{1}{q}} ) &=
\dfrac{ e^{t[1-x]_{\frac{1}{q}}} \, - \, e^{-t[1-x]_{\frac{1}{q}}}}{2} \\
&=  \dfrac{k!}{ 2(t[x]_q)^k} \left(  \dfrac{(t[x]_q)^k
e^{t[1-x]_{\frac{1}{q}}} - (t[x]_q)^k e^{-t[1-x]_{\frac{1}{q}} }}{k!}
\right)\\
&= \dfrac{k!}{ 2 t^k [x]_q^k} \left(  \sum_{n=k}^\infty B_{k, n}(
x \mid q) \dfrac{t^n}{n!} -  \sum_{n=k}^\infty B_{k, n}( x \mid
q)(-1)^{n-k} \dfrac{t^n}{n!} \right)\\
&=\dfrac{k!}{ 2 t^k [x]_q^k} \left(  \sum_{n=k}^\infty
(1-(-1)^{n+k} )  B_{k, n}( x \mid q) \dfrac{t^n}{n!} \right) \\
&=\dfrac{k!}{ 2  [x]_q^k} \left(  \sum_{n=0}^\infty (1-(-1)^{n} )
B_{k, n+k}( x \mid q) \dfrac{t^n}{(n+k)!} \right)  \\
&=\dfrac{1}{ 2  [x]_q^k} \left(  \sum_{n=0}^\infty (1-(-1)^{n} )
\dfrac{B_{k, n+k}( x \mid q)}{ \binom {n+k}{k}} \dfrac{t^n}{n!}
\right).\endaligned \eqno(12)
$$
Using Taylor expansion, we get

$$\dfrac{ e^{[1-x]_{\frac{1}{q}}t }-e^{-[1-x]_{\frac{1}{q}}t}}{2}
=\sum_{n=0}^\infty  \dfrac{[1-x]_{\frac{1}{q}}^n -(-[1-x]_{\frac{1}{q}})^n
}{2} \dfrac{t^n}{n!}. \eqno(13)
$$
Comparing the coefficients on the both sides of (12) and (13),
we obtain the following theorem.

\bigskip

{ \bf Theorem 6.}  For $   n, k  \in \Bbb Z_+$ and $ x \in [0,1]$,
we have
$$ [x]_q^{n+k}( [1-x]_{\frac{1}{q}}^n -(-1)^n [1-x]_{\frac{1}{q}}^n)
= [x]_q^n \big(1-(-1)^{n} \big) \dfrac{B_{k, n+k}( x \mid q)}{ \binom
{n+k}{k}}.$$

\bigskip

From (4) and Theorem 6, we note that
$$ \aligned  & B_{n, 2n+k}(1-x \mid \frac{1}{q})-(-1)^n B_{n, 2n+k}(1-x \mid
\frac{1}{q})\\
&=B_{n+k, 2n+k}(x \mid q)-(-1)^nB_{n+k, 2n+k}(x \mid q).
\endaligned \eqno(14)
$$
Therefore, by (14), we obtain the following theorem.

\bigskip

{ \bf Theorem 7.}  For $   n, k  \in \Bbb Z_+$ and $ x \in [0,1]$,
we have
$$ B_{n+k, 2n+k}(x \mid q)
= [x]_q^n \dfrac{B_{k, n+k}( x \mid q)}{ \binom {n+k}{n}}.$$

\bigskip
By Theorem 7, we get
$$ B_{n+1+k, 4n+2+k}(x \mid q)
=  \dfrac{B_{k, 2n+1+k}( x \mid q)}{ \binom
{2n+1+k}{2n+1}}[x]_q^{2n+1}.$$

 Thus, we have

$$\binom{2n+1+k}{2n+1}\binom{4n+2+k}{2n+1+k}=\binom{2n+1+k}{k}. $$

Let $ t \in \Bbb R$ and $ 0 \leq x \leq 1$. Then we see that
$$ \aligned   \cosh ( t[1-x]_{\frac{1}{q}} ) &=
\dfrac{ e^{t[1-x]_{\frac{1}{q}}}+e^{-t[1-x]_{\frac{1}{q}}}}{2} \\
&=  \dfrac{k!}{ 2(t[x]_q)^k} \left(  \dfrac{(t[x]_q)^k}{k!}e^{t[1-x]_{\frac{1}{q}}} + \dfrac{(-1)^k(-t[x]_q)^k
}{k!}e^{-t[1-x]_{\frac{1}{q}}}
\right)\\
&= \dfrac{k!}{ 2(t[x]_q)^k} \left(  \sum_{n=k}^\infty B_{k, n}( x
\mid q) \dfrac{t^n}{n!} +  \sum_{n=k}^\infty B_{k, n}( x \mid
q)(-1)^{n-k} \dfrac{t^n}{n!} \right)\\
&=\dfrac{1}{ 2  [x]_q^k} \left(  \sum_{n=0}^\infty (1+(-1)^{n} )
\dfrac{B_{k, n+k}( x \mid q)}{ \binom {n+k}{k}} \dfrac{t^n}{n!}
\right).\endaligned \eqno(15)
$$
Using Taylor expansion, we get

$$ \cosh( t[1-x]_{\frac{1}{q}} )
= \dfrac{1}{2} \sum_{n=0}^\infty (1+(-1)^n) [1-x]_{ \frac{1}{q}}^n
\dfrac{t^n}{n!}. \eqno(16)
$$
Therefore, by (15) and (16), we obtain
$$ [x]_q^{n+k}[1-x]_{\frac{1}{q}}^n = [x]_q^n  \dfrac{B_{k, n+k}( x \mid q)}{ \binom {n+k}{k}}.$$

Thus, we have

$$  B_{2n+k, 4n+k}(x \mid q)= \dfrac{B_{k, 2n+k}( x \mid q)}{ \binom {2n+k}{2n}}[x]_q^{2n}. \eqno(17)
$$
From (17), we have

$$\binom{4n+k}{2n+k}\binom{2n+k}{2n}=\binom{2n+k}{k}. $$

 \end{document}